\documentclass{article}
\usepackage{amsmath,amssymb,amsthm,enumerate}
\usepackage{graphics,graphicx,float,subcaption,caption,subcaption}
\usepackage{epstopdf}

\theoremstyle{definition}
\newtheorem{defn}{Definition} 

\begin{document}
%
%
 
\title{Solving fractional Hantavirus model: A new approach}
\author{Yogita Mahatekar \thanks{Department of Mathematics, College of Engineering Pune (COEP), Pune, India - 411005}   \thanks{Email: yogitasukale25@gmail.com}, Amey Deshpande  \thanks{School of Mathematics and Statistics, Dr. Vishwanath Karad MIT World Peace University, Pune, Maharashtra, India - 411038}  \thanks{2009asdeshpande@gmail.com} }

%
\maketitle

\begin{abstract}
A three equation differential biodiversity model depicting spread and propagation of Hantavirus epidemic amongst rodent species under the influence of alien competing mouse species was introduced by \cite{peixoto2006effect,hantavirus2010modeling}. We analyze the Caputo time-fractional version of the Biodiversity model along with added optimal control parameters as in \cite{mohamed2019modeling}. A new algorithm is proposed to solve fractional differential equations which is derived from the implementation of NIM \cite{daftardar2006iterative} and implicit $\theta$-method \cite{yakit2018explicit} with $\theta=1$. This newly proposed algorithm is then used to integrate the Optimal controlled Caputo time-fractional version of Biodiversity Hantavirus model (OC-CTFBH model). The obtained solutions are plotted and analyzed further for the effects of parameter variations. 
\end{abstract}
\section{Introduction}
In January 2020, when the World Health Organization (WHO) first mentioned a cluster of pneumonia cases in Wuhan, no one expected that by the $11^{th}$ of March, COVID-19 (Coronavirus Disease-19) would have spread across the globe and would be declared as an pandemic. The well-known works of W. O. Kermack and A. G. McKendrick, published in 1927, on the compartmental SIR model, have paved the way for epidemiologists to mathematically describe the dynamics of infectious diseases. These SIR models have become increasingly popular among the scientific community, particularly for the analysis of current pandemic. These compartmental SIR-models are comprised of systems of ordinary / fractional differential equations that, except for a few particular cases, do not admit analytical solutions. For the complete dynamical analysis and study of compartmental models, it is important to fully resolve corresponding differential equations / fractional differential equations. And therefore it is pertinent to develop novel numerical methods to obtain best possible approximation to the exact solution. These numerical methods and obtained solutions provide a fair analysis of the spread of the epidemic and enables us to make more robust projections on the basis of of past observations.

A three equation differential biodiversity model depicting spread and propagation of Hantavirus epidemic amongst rodent species under the influence of alien competing mouse species was introduced by \cite{peixoto2006effect,hantavirus2010modeling}. In this experiment the species of the rodent which are susceptible to Hantavirus epidemic are made to compete for the resources with the `alien' mouse species externally introduced in the same environment. The proposed model is termed as Biodiversity model and authors \cite{peixoto2006effect} analyze and conclude that the increased competitive pressure will likely result into reduction or complete elimination of the prevalence of infection among the rodent species. An optimal control was introduced to this model by \cite{mohamed2019modeling} where a new parameter $E$ was introduced depicting `Harvesting efforts' which amount to removal of rodents periodically which happens naturally due to various factors such as hunting or capturing. This new optimal control parameter increases chance of the system to stabilize towards \emph{co-existence equilibrium point} and indeed helps in reduction of the infection among rodents. In this article, we analyze dynamics of Optimal controlled Caputo time-fractional version of Biodiversity Hantavirus model (OC-CTFBH model).  To do this, we introduce a novel numerical method developed with a combination of implicit form of $\theta-$method and a decomposition technique developed by Daftardar-Gejji and H. Jafari namely NIM \cite{daftardar2006iterative}. The obtained solutions by this newly proposed method are compared with the exact and/or numerical solutions obtained by other prevalent methods. Figures representing susceptible rodents, infected rodents, and alien populations are plotted. 
It can be observed that susceptible rodent population initially  increases sharply, then reaches to a maximum value and later approaches to an equilibrium quantity. Similar behavior is observed in case of infected and alien populations. It has been observed that, when we apply optimal control namely harvesting efforts to the model then both susceptible and infected rodent population minimizes. 

The article is further organized as follows. Section \ref{preliminaries} introduces to basic terminology and definitions from fractional calculus used in this article. Section \ref{formulation} discuss various Hantavirus models leading to formulation of OC-CTFBH model. Consequent Sections \ref{numericalmethod} and Section \ref{numericalsimulation} proposes new numerical method after taking a brief overview of NIM and implicit $\theta$-method and numerical simulations of OC-CTFBH using this method respectively. Section \ref{conclusions} sums up the findings and observations of the article.

\section{Preliminaries} \label{preliminaries}
In this section, we introduce some preliminaries from fractional calculus. For more details, we refer the readers to  \cite{machado2011recent, diethelm2010analysis}.
\begin{defn}
	Grunwald Letnikov Fractional derivative operator $_{GL}D^{\alpha}$ is defined as  
	\begin{equation}
		_{GL}D^{\alpha} f(t)= \lim\limits_{N \rightarrow \infty} h^{-\alpha}_{N}\sum^{N}_{j=0} w^{(\alpha)}_{j} f(t-j_N);\, ~\text{ where }~
		w^{(\alpha)}_{j} = \frac{\Gamma(j-\alpha)}{\Gamma(-\alpha)\Gamma(j+\alpha)} 
	\end{equation} 
\end{defn}

\begin{defn}
	The Riemann-Liouville fractional integral of order $\alpha > 0$ of $f \in C[0, \infty)$ is defined as
	\begin{align}
	I^{\alpha}f(t)=\frac{1}{\Gamma(\alpha)}\int\limits_{0}^{t}\frac{f(\tau)}{(t-\tau)^{1-\alpha}}d\tau.
	\end{align}
\end{defn}
\begin{defn}\cite{diethelm2010analysis}\label{rlder}
	Let $\alpha \in \mathbb{R^{+}}$ and $m = \lceil{\alpha}\rceil.$ The Riemann-Liouville fractional operator of order $\alpha$, $_{RL}D^{\alpha}$, is defined as
	  \begin{equation}
	  _{RL}D^{\alpha}f(t) = \frac{d^m}{dt^m}I^{m-\alpha}f(t).
\end{equation}
\end{defn}
\noindent Note that for $\alpha \in \mathbb{R^{+}}$ and $p \in \mathbb{N}$ such that $p > \alpha.$ Then $_{RL}D^{\alpha}f(t)=\frac{d^p}{dt^p}\;I^{p-\alpha}f(t)$.
\begin{defn}
	The Caputo derivative of order $\alpha \in (k-1, k], ~ k \in \mathbb{N}$ of $f \in C^k(0, \infty)$ is defined as 
	\begin{equation}
	^{c}D^{\alpha}f(t) =
	\begin{cases}
	
	\frac{1}{\Gamma(k-\alpha)} ~\int_{0}^{t} (t-\tau)^{k-\alpha-1} f^{(k)}(\tau)\: d\tau,~~ &\alpha \in (k-1,k),  \\
	 f^{(k)}(t), ~~ & \alpha = k.
	 \end{cases}
	\end{equation}
\end{defn}
\noindent The relation between RL and Caputo fractional derivative is as follows. 
\begin{equation}
	^{c}D^{\alpha} f(t)= _{RL}D^{\alpha} (f(t)-f(t_0))\label{eq10}
\end{equation}
and under suitable regularity assumption, 
\begin{equation}
	_{RL}D^{\alpha} f(t)= 	_{GL}D^{\alpha} f(t)\label{eq11}
\end{equation}
And hence in view of equations (\ref{eq10}) and (\ref{eq11}) we have
\begin{equation}
	^{c}D^{\alpha} f(t)= _{RL}D^{\alpha} (f(t)-f(t_0))=_{GL}D^{\alpha} (f(t)-f(t_0))\label{eq12}
\end{equation}
Therefore; 
\begin{equation}
	^{c}D^{\alpha} f(t)=_{GL}D^{\alpha}(f(t)-f(t_0))=h^{-\alpha}\sum^{n}_{j=0}w^{(\alpha)}_{j}(f(t-j)-f(t_0))\label{eq133}
\end{equation}
Hereafter in the article, to simplify the notations, we drop the left superfix $c$ and use notation $D^{\alpha}f(t)$ to denote Caputo fractional derivative operator. 
%
\section{Formulation of Hanta virus model} \label{formulation}
The Biodiversity model was proposed by Peixoto and Abramson \cite{peixoto2006effect}. In this model population of the rodent is introduced with Hantavirus infection and is allowed to compete for the resources with an alien uninfected rodent species. 

Let the variables $x(t), y(t), z(t)$ denote the population of rodent susceptible to Hantavirus, population infected with Hantavirus and the alien rodent population at time $t$ respectively. The $r(t) = x(t) + y(t)$ represents the total population of the rodents introduced with Hantavirus. Then the Biodiversity model is given as follows. 
  
\begin{align}
	\frac{dx}{dt}&=br-cx-\frac{xr}{k(t)}-axy \\
	\frac{dy}{dt}&=-cy-\frac{yr}{k(t)}+axy
\end{align}
where
$b:$ birth rate,\,$c:$ natural mortality rate,\,$a:$ transmission rate(agression parameter),\,$k:$ environmental parameter.

 Yusof et.al.\cite{mohamed2019modeling} introduces optimal control parameter $E$ representing `Harvesting efforts' to Biodiversity model. This parameter models periodical reduction of population due to various natural instances such as hunting, capturing into the model and optimal control helps in stabilization of the system. This model is given as follows.    
\begin{align}
	\frac{dx}{dt}&=br-cx-\frac{x(r+qz)}{k(t)}-axy-E x \\
	\frac{dy}{dt}&=-cy-\frac{y(r+qz)}{k(t)}+axy-Ey \\
	\frac{dz}{dt}&=(\beta-\gamma)z-\frac{z}{k(t)}(z+\epsilon r) 
\end{align} 
where, $E:$ Harvesting efforts,\,$q:$ influence of the alien population,\, $z(t):$ population of alien,\,$\beta, \gamma,\epsilon:$ corresponding parameters to obtain resources from other species,\,$k:$ environmental parameter.

The Fractional version of the corresponding model gives extra parameter of flexibility and control to the existing model. This allows us to more accurately depict the real life dynamics and match with the observed experimental data more precisely. 
%
Therefore for fractional order $0 < \alpha \leq 1$, the  Optimal controlled Caputo time-fractional version of Biodiversity Hantavirus model (OC-CTFBH model) is given as follows.
\begin{align} 
D^{\alpha}x &= b r - c x - \frac{x}{k}\; (r + q z) - a x y - E x \label{model11} \\
D^{\alpha}y &= -c y - \frac{y}{k}\; (r + q z) + a x y - E y \label{model12}\\
D^{\alpha}z &= (\beta - \gamma) z - \frac{z}{k} \;(z + \epsilon \;r) \label{model13} 
\end{align} where $D^{\alpha}f(t)$ represents Caputo fractional derivative of $f(t)$ of order $\alpha$.  The parameters and variables in OC-CTFBH represents as follows.
\begin{itemize}
	\item $x \equiv x(t) :$ the susceptible rodent population.
	\item $y \equiv y(t) :$ the infected rodent population. $r(t) = x(t) + y(t)$.
	\item $z \equiv z(t) :$ the alien population.
	\item $b, c : $ the birth and death rate of the rodent respectively.
	\item $\beta, \gamma : $ the birth and death rate of the alien respectively.
	\item $q : $ influence exerted by the alien population on the rodent population. Note that $q : 0$ will reduce the Biodiversity model to the model without a competition represented by first two equations.
	\item $\epsilon : $ the influence exerted by the rodent population on the alien species.
	\item $a : $ propagation rate of the Hantavirus. 
	 \item $k : $ represents the carrying capacity of the environment. 
\end{itemize}   

\section{Numerical Algorithm} \label{numericalmethod}
In this section we propose a new numerical method as an combination of NIM and implicit $\theta$-method.  
\subsection{New iterative method (NIM)}\label{nim}
Daftardar-Gejji and Jafari \cite{daftardar2006iterative} developed a new iterative method to solve functional equations of the form:
\begin{equation}
	v=h + L(v)+ N(v)\label{em7}
\end{equation}
where $h$ is a known part, $L$ is a linear operator while $N$ a non linear operator.

In this method, it is assumed that equation \eqref{em7} has a solution of the form 
$v=\sum^{\infty}_{i=0}v_{i}$ where 
 \begin{align*}
 	v_0 &= h,\\
 	v_1 &= L(v_0)+N(v_0), \\
 	v_i &= L(v_{i-1}) + N\left( \sum_{j=0}^{i-1}\; v_j\right) - N\left( \sum_{j=0}^{i-2}\; v_j\right), ~~~i=2,3,\cdots
 \end{align*}
Since $L$ is a linear operator, we have $L(v)=L\left(\sum^{\infty}_{i=0}v_{i}\right)=\sum^{\infty}_{i=0}L(v_{i})$ and non linear operator $N$ is decomposed as
\begin{equation*}
N(v)=N(v_{0})+ \sum_{i=2}^{\infty} \left[ N\left( \sum_{j=0}^{i-1}\; v_j\right) - N\left( \sum_{j=0}^{i-2}\; v_j\right)\right] 
\end{equation*}
Therefore we get
 \begin{align*}
 v&=v_{0}+v_{1}+v_{2}+\cdots \\&=
h+L(v_0)+N(v_{0})+L(v_1)+[N(v_{0}+v_{1})-N(v_{0})]+\cdots\\&=h+L(v)+N(v),
\end{align*}
 and hence $v$
satisfies the functional equation \eqref{em7}. A general $k$-term NIM solution is given by $v=\displaystyle\sum^{k-1}_{i=0}v_{i}.$

\subsection{Novel $\theta$-numerical method}
A general form of implicit $\theta$-method is given as 
\begin{equation}\label{eqqqqq}
	D^{\alpha} y(t_n)=\theta f(t_n, y_n)+(1-\theta) f(t_{n-1},y_{n-1}).
\end{equation}
In view of equation \eqref{eq133}, we have 
\begin{align} 
	h^{-\alpha}\sum^{n}_{j=0}w^{(\alpha)}_{j}(y(t_{n-j})-y(t_0)) &= \theta f(t_n, y_n)+(1-\theta) f(t_{n-1},y_{n-1})\\ 
	\implies h^{-\alpha}\sum^{n}_{j=0}w^{(\alpha)}_{j}(y_{n-j}-y_0) &= \theta f(t_n, y_n)+(1-\theta) f(t_{n-1},y_{n-1}) 
\end{align}
In particular, for $\theta=1$ we get
\begin{equation}
	h^{-\alpha}\sum^{n}_{j=0}w^{(\alpha)}_{j}(y_{n-j}-y_0) = f(t_n, y_n)\label{eq17}
\end{equation}
Applying the novel implicit $\theta$-method as given in \eqref{eqqqqq} to \eqref{eq17} to solve the OC-CTFBH model \eqref{model11}-\eqref{model13} we get
\begin{align} 
	h^{-\alpha}\sum^{n}_{j=0}w^{(\alpha)}_{j}(x_{n-j}-x_{0}) &= (b-c)x_{n}+by_{n}-\frac{x_{n}}{k}(x_{n}+y_{n}+qz_{n})\label{eq18}\\
	\nonumber &~~~~~~~~~~~~~~~~~~~~~~~~~~~~~~~~~~~~~~~~~-ax_{n}y_{n}-Ex_{n} \\
	h^{-\alpha}\sum^{n}_{j=0}w^{(\alpha)}_{j}(y_{n-j}-y_{0}) &= (-c)y_{n}-\frac{y_{n}}{k}(x_{n}+y_{n}+qz_{n})\label{eq19}\\
	\nonumber &~~~~~~~~~~~~~~~~~~~~~~~~~~~~~~+ax_{n}y_{n}-Ey_{n} \\
	h^{-\alpha}\sum^{n}_{j=0}w^{(\alpha)}_{j}(z_{n-j}-z_{0}) &= (\beta-\gamma)z_{n}-\frac{z_{n}}{k}(z_{n}+\epsilon x_{n}+\epsilon y_{n})\label{eq20}
\end{align}
Equation (\ref{eq18}) we express in the form of NIM namely $u=N(u)+f$ where $f$ is the known part, $N$ is the non-linear operator and $u$ is the unknown. Three term NIM solution of this functional equation is given by $u=u_0+u_1+u_2$ where $u_0=f,\,u_1=N(u_0),\,u_{2}=N(u_0+u_1)-N(u_0)$. Now expressing (\ref{eq18}) in the form $u=N(u)+f$ as follow
\begin{equation}
	\begin{split} 
	h^{-\alpha}\sum^{n}_{j=1}w^{(\alpha)}_{j}(x_{n-j}-x_{0})+h^{-\alpha}w^{(\alpha)}_{0}(x_{n}-x_{0}) = (b-c)x_{n}+by_{n}-\frac{r^{2}_{s_{n}}}{k}\\ -\frac{x_{n}y_{n}}{k}-\frac{qx_{n}z_{n}}{k}-ax_{n}y_{n}-Ex_{n}\label{eq21}
\end{split}
\end{equation}
and upon rearranging terms, we get 
\begin{equation}
\begin{split} 
	(h^{-\alpha}w^{(\alpha)}_{0}-b+c+E)x_{n}=	-h^{-\alpha}\sum^{n}_{j=1}w^{(\alpha)}_{j}(x_{n-j}-x_{0})+h^{-\alpha}w^{(\alpha)}_{0}x_{0}\\ + by_{n}-\frac{r^{2}_{s_{n}}}{k}-\frac{x_{n}y_{n}}{k}-\frac{qx_{n}z_{n}}{k}-ax_{n}y_{n}\label{eq22}
\end{split}
\end{equation}
Noting that equation \eqref{eq22}
is of the form $x_{n}=N_1(x_{n},y_{n},z_{n})+f_{1_{n}}$ where
\begin{equation}
	f_{1_{n}}=\frac{-h^{-\alpha}\sum^{n}_{j=1}w^{(\alpha)}_{j}(x_{n-j}-x_{0})+h^{-\alpha}w^{(\alpha)}_{0}x_{0}}{h^{-\alpha}w^{(\alpha)}_{0}-b+c+E}\label{eq23}
\end{equation}
and \begin{equation}
	N_1(x_{n},y_{n},z_{n})=\frac{by_{n}-\frac{r^{2}_{s_{n}}}{k}-\frac{x_{n}y_{n}}{k}-\frac{qx_{n}z_{n}}{k}-ax_{n}y_{n}}{h^{-\alpha}w^{(\alpha)}_{0}-b+c+E}\label{eq24}
\end{equation}
Similarly we arrange equation \eqref{eq19} in the form say $y_{n}=N_{2}(x_{n},y_{n},z_{n})+f_{2_{n}}$ as follows:
\begin{equation}
	\begin{split} 
	h^{-\alpha}\sum^{n}_{j=1}w^{(\alpha)}_{j}(y_{n-j}-y_{0})+h^{-\alpha}w^{(\alpha)}_{0}(y_{n}-y_{0}) = (-c)y_{n}-\frac{r^{2}_{i_{n}}}{k}-\frac{y_{n}x_{n}}{k}\\-\frac{qy_{n}z_{n}}{k}+ax_{n}y_{n}-Ey_{n}\label{eq25}
\end{split}
\end{equation}
which implies that:
\begin{equation}
	\begin{split} 
	(h^{-\alpha}w^{(\alpha)}_{0}+c+E)y_{n}=	-h^{-\alpha}\sum^{n}_{j=1}w^{(\alpha)}_{j}(y_{n-j}-y_{0})+h^{-\alpha}w^{(\alpha)}_{0}y_{0}-\frac{r^{2}_{i_{n}}}{k}\\-\frac{x_{n}y_{n}}{k}-\frac{qy_{n}z_{n}}{k}+ax_{n}y_{n}\label{eq26}
\end{split}
\end{equation}
Here we denote 
\begin{equation}
	f_{2_{n}}=\frac{-h^{-\alpha}\sum^{n}_{j=1}w^{(\alpha)}_{j}(y_{n-j}-y_{0})+h^{-\alpha}w^{(\alpha)}_{0}y_{0}}{h^{-\alpha}w^{(\alpha)}_{0}+c+E}\label{eq27}
\end{equation}
and \begin{equation}
	N_{2}(x_{n},y_{n},z_{n})=\frac{-\frac{r^{2}_{i_{n}}}{k}-\frac{x_{n}y_{n}}{k}-\frac{qy_{n}z_{n}}{k}+ax_{n}y_{n}}{h^{-\alpha}w^{(\alpha)}_{0}+c+E}\label{eq28}
\end{equation}

\noindent In a similar manner, equation \eqref{eq20} is arranged in the form of $z_{n}=N_{3}(x_{n},y_{n},z_{n})+f_{3_{n}}$ as follows:
\begin{equation}
	\begin{split} 
	h^{-\alpha}\sum^{n}_{j=1}w^{(\alpha)}_{j}(z_{n-j}-z_{0})+h^{-\alpha}w^{(\alpha)}_{0}(z_{n}-z_{0}) = (\beta-\gamma)z_{n}-\frac{z^{2}_{a_{n}}}{k}\\-\frac{\epsilon z_{n}x_{n}}{k}-\frac{z_{n}\epsilon y_{n}}{k}\label{eq29}
	\end{split}
\end{equation}
and here we denote
\begin{equation}
	f_{3_{n}}=\frac{-h^{-\alpha}\sum^{n}_{j=1}w^{(\alpha)}_{j}(z_{n-j}-z_{0})+h^{-\alpha}w^{(\alpha)}_{0}z_{0}}{h^{-\alpha}w^{(\alpha)}_{0}-\beta+\gamma}\label{eq30}
\end{equation}
and \begin{equation}
	N_{3}(x_{n},y_{n},z_{n})=\frac{-\frac{z^{2}_{a_{n}}}{k}-\frac{\epsilon z_{n}x_{n}}{k}-\frac{z_{n}\epsilon y_{n}}{k}}{h^{-\alpha}w^{(\alpha)}_{0}-\beta+\gamma}\label{eq31}
\end{equation}
\subsection{Finding solution $x(t),\,y(t),\,z(t)$ using NIM} 
\begin{align}
	x_{n}&=x(t_{n})=\sum_{k=0}^{2} x_{kn}=x_{0n}+x_{1n}+x_{2n}\label{eq32}\\
	y_{n}&=y(t_{n})=\sum_{k=0}^{2} y_{kn}=y_{0n}+y_{1n}+y_{2n}\label{eq33}\\
	z_{n}&=z(t_{n})=\sum_{k=0}^{2} z_{kn}=z_{0n}+z_{1n}+z_{2n}\label{34}	
\end{align}

\noindent Note that when $n=0,$ we are at time parameter $t=t_0$ and in this case $x_{0}=x(t_{0}),\,y_{0}=y(t_{0}),\,z_{0}=z(t_{0})$ which are the given initial conditions in the model.
When $n=1$, we are finding $x(t),\,y(t),\,z(t)$ at time $t=t_{1}.$ Using 3-term NIM solution, $x_{1}=x(t_{1}),\,y_{1}=y(t_{1}),\,z_{1}=z(t_{1})$  are obtained as follows:
\begin{align}
	x_{1}&=x(t_{1})=\sum_{k=0}^{2} x_{k1}=x_{01}+x_{11}+x_{21}\label{eq35}\\
	y_{1}&=y(t_{1})=\sum_{k=0}^{2} y_{k1}=y_{01}+y_{11}+y_{21}\\
	z_{1}&=z(t_{1})=\sum_{k=0}^{2} z_{k1}=z_{01}+z_{11}+z_{21}\label{37}	
\end{align}
where $x_{01}=f_{11};\,y_{01}=f_{21};\,z_{01}=f_{31}.$ Further 
$x_{11}=N_{1}(f_{11},f_{21},f_{31});\,y_{11}=N_{2}(f_{11},f_{21},f_{31});\,z_{11}=N_{3}(f_{11},f_{21},f_{31})$ and 
\begin{align}
x_{21}&=N_{1}(x_{01}+x_{11},y_{01}+y_{11},z_{01}+z_{11})-N_{1}(f_{11},f_{21},f_{31})\\
y_{21}&=N_{2}(x_{01}+x_{11},y_{01}+y_{11},z_{01}+z_{11})-N_{2}(f_{11},f_{21},f_{31})\\
z_{21}&=N_{3}(x_{01}+x_{11},y_{01}+y_{11},z_{01}+z_{11})-N_{3}(f_{11},f_{21},f_{31})
\end{align}
\section{Numerical simulations}\label{numericalsimulation}
In this section, we discuss numerical results obtained by solving Hantavirus model using new iterative method and implicit $\theta-$ method $(\theta=1)$. Parameter values used in the numerical simulations are: $b=1,\,c=0.6,\,a=0.1,\,q=0.2,\,\beta=1,\,\gamma=0.5,\,\epsilon=0.1.$ $k$ is varied which represents the environmental factor. Here system of equations representing Hantavirus model is solved for $k=20$ and for $k=250$ with initial conditions $x(0)=5,\,y(0)=5,\,z(0)=5.$  
Step length $h=0.01.$
\begin{figure}[H]
	\centerline{\includegraphics[width=250pt]{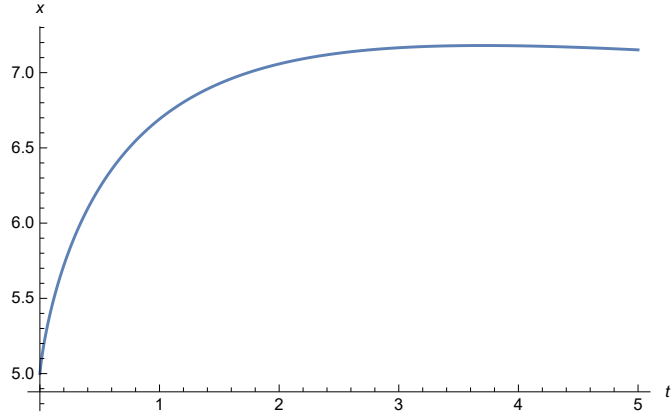}}
	\caption{ x(t) versus t (in months) for $k=20,\,E=0$ }\label{fig1}		
\end{figure}
\begin{figure}[H]
	\centerline{\includegraphics[width=250pt]{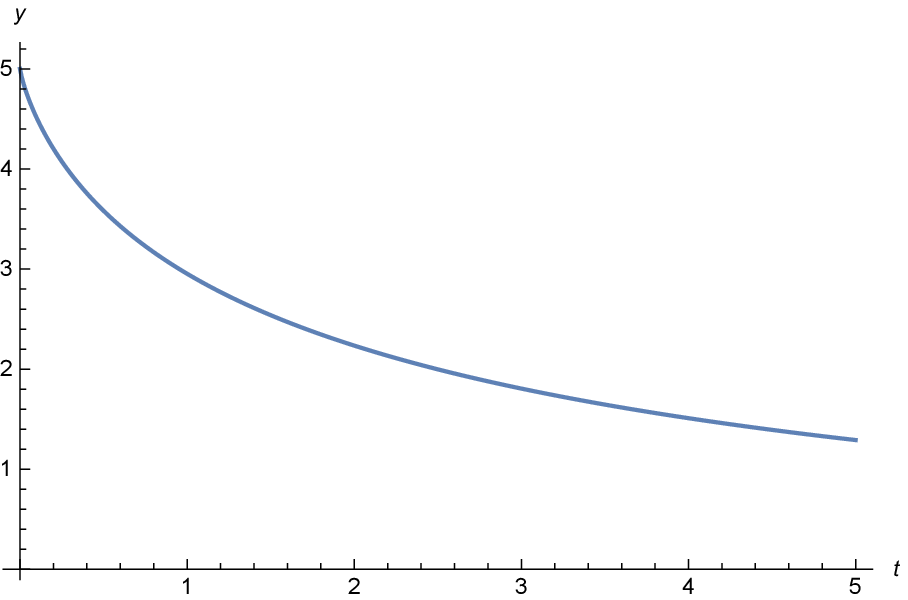}}
	\caption{y(t) versus t (in months) for $k=20,\,E=0$ }\label{fig2}		
\end{figure}
\begin{figure}[H]
	\centerline{\includegraphics[width=250pt]{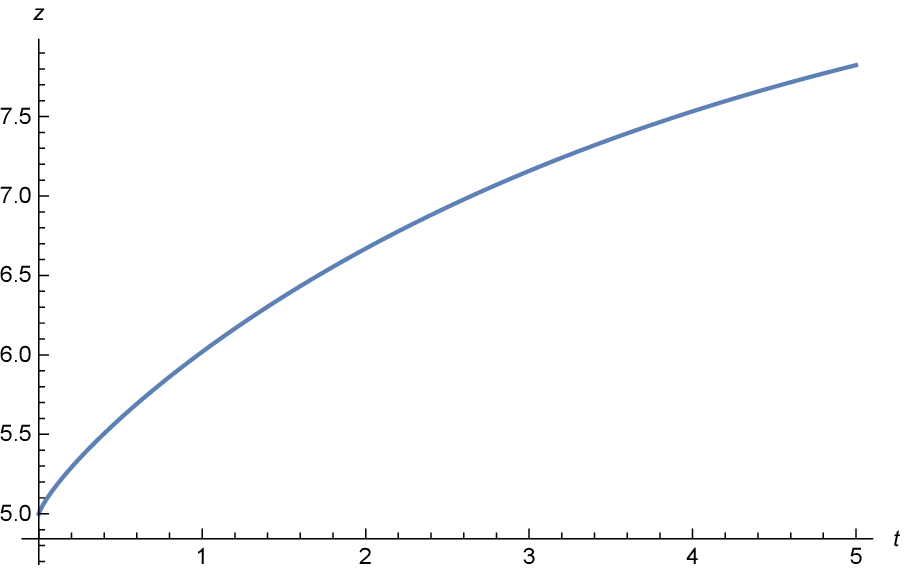}}
	\caption{z(t) versus t (in months) for $k=20,\,E=0$ }\label{fig3}		
\end{figure}
\begin{figure}[H]
	\centerline{\includegraphics[width=250pt]{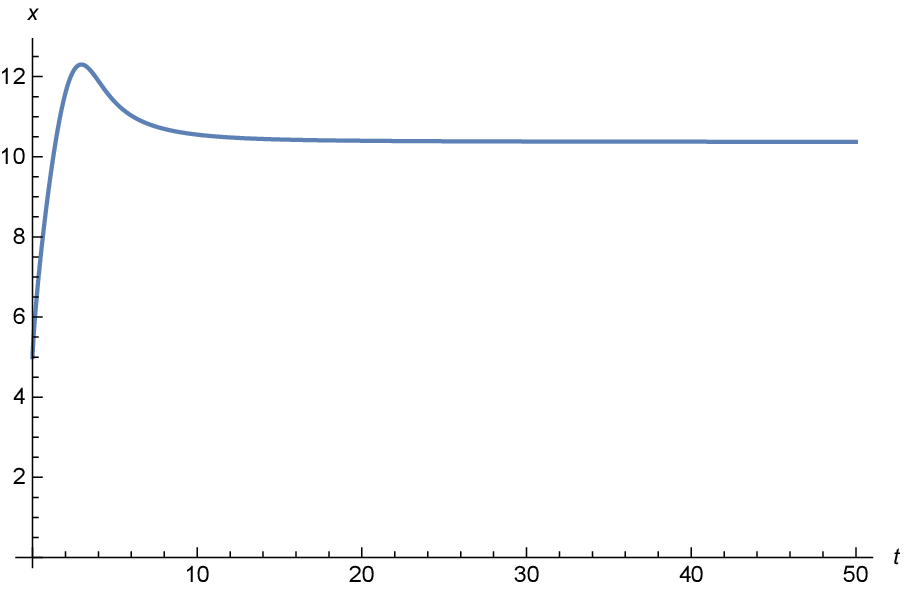}}
	\caption{x(t) versus t (in months) for $k=250,\,E=0$ }\label{fig4}		
\end{figure}
\begin{figure}[H]
	\centerline{\includegraphics[width=250pt]{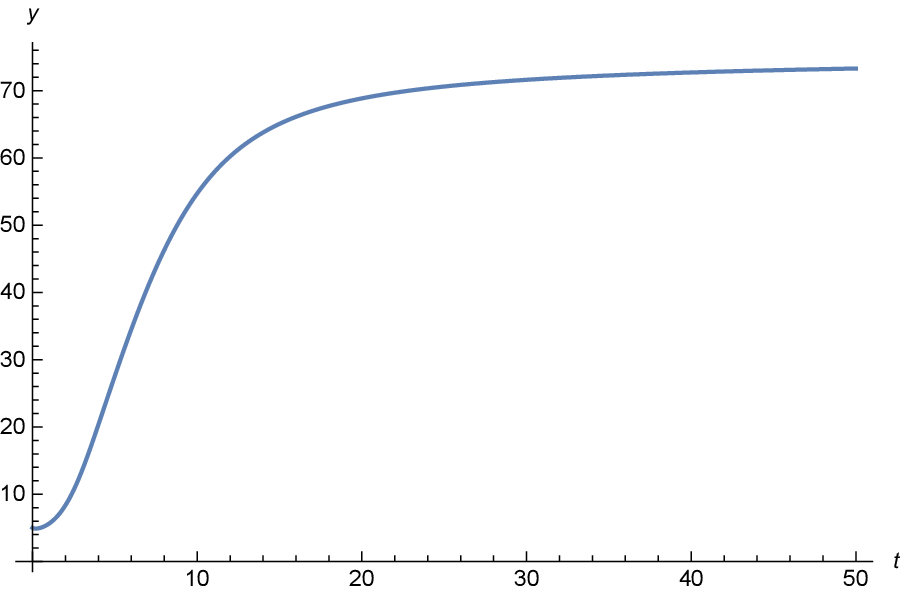}}
	\caption{y(t) versus t (in months) for $k=250,\,E=0$ }\label{fig5}		
\end{figure}
\begin{figure}[H]
	\centerline{\includegraphics[width=250pt]{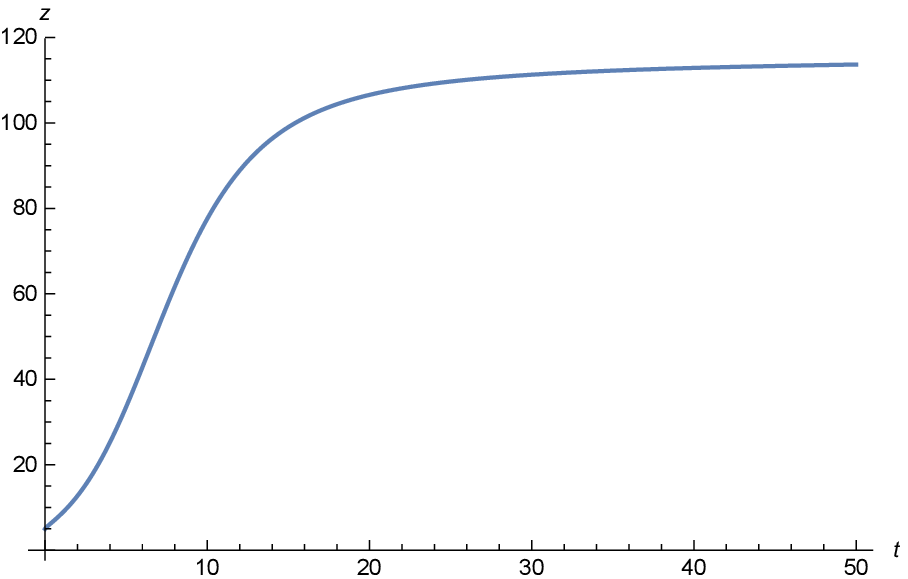}}
	\caption{z(t) versus t (in months) for $k=250,\,E=0$ }\label{fig6}		
\end{figure}

\begin{figure}[H]
	\centerline{\includegraphics[width=250pt]{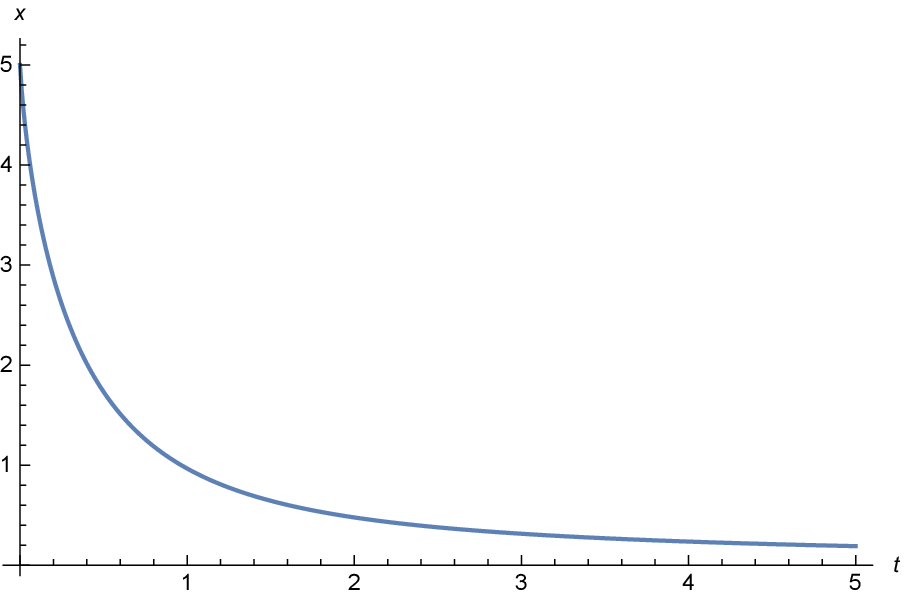}}
	\caption{x(t) versus t (in months) for $k=20,\,E=0.5$}	\label{fig7}	
\end{figure}
\begin{figure}[H]
	\centerline{\includegraphics[width=250pt]{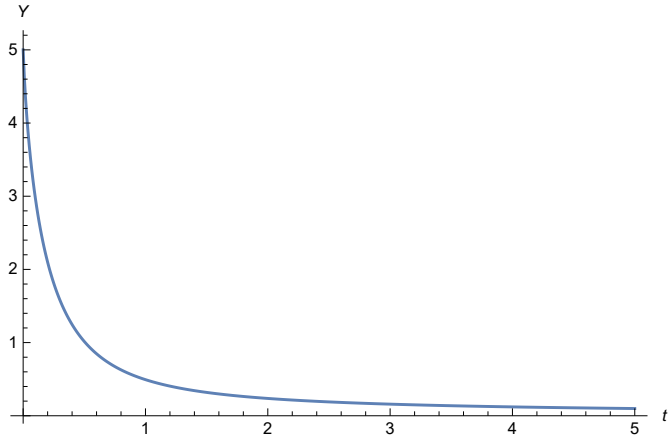}}
	\caption{ y(t) versus t (in months) for $k=20,\,E=0.5$}	\label{fig8}	
\end{figure}

The dynamics of the Hantavirus model without optimal control (taking the harvesting efforts $E=0$)  and with optimal control by taking the harvesting efforts $E=0.5$ are shown in Figs.(\ref{fig1}-\ref{fig8}).
In Figs.\ref{fig1}, \ref{fig2}, \ref{fig3},\ref{fig4},\ref{fig5},\ref{fig6} the population of susceptible rodent, infected rodent and alien population is depicted without optimal control ($E=0$) for for various values of $k.$ Wheareas in Fig.\ref{fig7} and in Fig. \ref{fig8}  the susceptible rodent population $x(t)$ and infected rodent population $y(t)$ is shown with optimal control by taking harvesting effects parameter $E=0.5.$ It has been observed that the susceptible rodent population $x(t)$ and infected rodent population $y(t)$ will disappear and stabilize at a certain steady value, when the control parameter $E$ is applied. The susceptible rodent population and even infected population is observed to be behaving oppositely in the model without optimal control. These susceptible and infected populations are initially increasing sharply and then reaches to a certain maximum value before approaching to an equilibrium value.

\section{Conclusions}\label{conclusions}
In this paper, we developed a new algorithm to solve fractional order Hantavirus model with and without optimal control. New numerical algorithm is accurate and efficient to apply for solving fractional order differential equations representing Hantavirus model. New algorithm is a combination of implicit $\theta-$ method and new iterative method. We studied rich dynamics of  Hantavirus model using this new numerical technique and effects of harvesting effects $E(t)$ as an optimal control on spread of Hantavirus infection is studied. The theoretical study and numerical simulations of the model are discussed in this paper. It has been observed that, when we apply optimal control namely harvesting efforts to the model then both susceptible and infected rodent population minimizes. This clearly states that optimal control is able to eliminate infected rodent population with the presence of harvesting effects.
\bibliographystyle{abbrv}
\bibliography{paper-arxive}
\end{document}